\documentclass[ukrainian,12pt]{article}
\usepackage{amsmath}
\usepackage{amsbsy}
\usepackage{amssymb}
\usepackage{amscd}
\usepackage{a4wide}
\usepackage[T2A]{fontenc}
\usepackage[cp1251]{inputenc}
\usepackage[ukrainian,english]{babel}
\usepackage{epic}

\newtheorem{theorem}{Theorem}

\newtheorem{definition}{Definition}

 \DeclareMathOperator{\rep}{Rep}
 \DeclareMathOperator{\image}{\Im m}
 \DeclareMathOperator{\ob}{Ob}

\begin{document}

\title{ The spectral  problem and *-representations of algebras associated with Dynkin graphs.}
\author{\footnotesize STANISLAV KRUGLJAK, STANISLAV POPOVYCH, YURII SAMOILENKO \footnote{Insitute of Mathematics of National Academy of Sciences, Kyiv, Ukraine}}

\maketitle

\newcommand{\ba}[1]{\ensuremath{e_{\phi(#1)}}}
\newcommand{\sca}[2]{\ensuremath{\langle #1,#2 \rangle}}
\newcommand{\algebra}[1]{\ensuremath{\mathbb{C}\langle #1 \rangle}}
\newcommand{\algp}[1]{\ensuremath{\mathcal{P}_{n,#1}}}
\newcommand{\pabo}[1]{\ensuremath{\mathcal{P}_{n,abo,#1}}}
\newcommand{\g}[2]{\ensuremath{\gamma_{#1}^{(#2)}}}
\newcommand{\ma}[2]{\ensuremath{\alpha_{#1}^{(#2)}}}
\numberwithin{equation}{section}
 \abstract{
We study the connection between *-representations of algebras associated with
graphs, locally-scalar graph representations and the problem about the  spectrum of a  sum of two Hermitian operators. For algebras associated with Dynkin graphs we give an  explicit description of the  parameters for which   there are  irreducible representations and an algorithm for contructing  these  representations. 


 \medskip\par\noindent
AMS SUBJECT CLASSIFICATION: 16W10, 16G20, 47L30 
  }
\vspace{10mm}

\noindent {\Large \bf Introduction.} \vspace{5mm}

 {
 1. Let  $A$, $B$, $C$  be Hermitian  $n\times n$ matrices with given
eigenvalues, $\tau(A)= \{\lambda_1(A )\ge \lambda_2(A)\ge \ldots
\ge \lambda_n(A) \}$, $\tau(B)= \{\lambda_1(B )\ge \lambda_2(B)\ge
\ldots \ge \lambda_n(B)\}$, $\tau(C)= \{\lambda_1(C )\ge
\lambda_2(C)\ge \ldots \ge \lambda_n(C) \}$. The well-known
classical problem  about spectrum of a sum of two Hermitian
matrices (Horn's problem) is to describe a connection between
 $\tau(A), \tau(B), \tau(C)$ for matrices $A$, $B$, $C$ such that
$A+B=C$.

A recent solution of this problem  (see. ~\cite{fulton, klychko}
 and others) gives a complete description of
possible $\tau(A)$, $\tau(B)$, $\tau(C)$ in terms of linear
inequalities of the form
 \begin{equation}\label{eq1}
 \sum_{i\in I}\lambda_i(A)
+\sum_{j\in J}\lambda_j(B )\ge \sum_{k\in K} \lambda_k(C),
\end{equation}
where  $I, J, K$ are certain subsets of  $\{1,\ldots, n\}$.
 Note that the number of necessary inequalities increases with
$n$.

\noindent 2. Consider the following modifications of the problem
mentioned above,  called henceforth the {\it spectral 
problem} (res. {\it strict spectral  problem}). We will consider bounded linear Hermitian operators on a separable Hilbert space  $H$. For an operator   $X$ denote by $\sigma(X)$ its spectrum. Let  $M_1, M_2, M_3$ be given closed
subsets of $\mathbb{R}^+$ and $\gamma\in\mathbb{R}^+$. The problem
consists of the following:  1)to determine whether there are Hermitian operators $A, B, C$ on  $H$ such that $\sigma(A)\subseteq M_1$, $\sigma(B)\subseteq M_2$, $\sigma(C)\subseteq M_3$(res. $\sigma(A)= M_1$, $\sigma(B)= M_2$, $\sigma(C)= M_3$) and $A+B+C=\gamma I$?  2) if  the answer is in the  affirmative, to  give a description (up to unitary equivalence) of the operators. In this work the sets $M_1, M_2, M_3$ will be finite. Note that even for finite $M_k$ the second part of the problem can be very complicated if  $|M_k|$ is large enough. 

The essential difference with the classical Horn's problem is that we
do not fix the dimension of  $H$   (it may be finite or infinite)
  and the spectral multiplicities. It seems that  solution of spectral and strict spectral  problems could not be deduced directly from inequalities of the form (\ref{eq1}), since the number of necessary inequalities increases with $n$. 

\noindent 3. These problems could be reformulated in terms of *-representations of *-algebras.  Namely, let $\alpha=(\alpha_1,\alpha_2,\ldots, \alpha_{k})$,   $\beta=(\beta_1, \beta_2,\ldots, \beta_{l})$,   $\delta=(\delta_1, \delta_2,\ldots,\delta_{m})$ be vectors with positive   strictly decreasing coefficients. Let us remark that we can assume that each set  $M_1, M_2, M_3$ contains zero (this can be achieved by a translation). Henceforth,   $M_1=\alpha\cup\{0\}$, $M_2=\beta\cup\{0\}$, $M_3=\delta\cup\{0\}$. Let us consider the  associative algebra defined by
the following generators and relations  (see.~\cite{popsam}):
\begin{gather*}
\mathcal{P}_{\alpha,\beta,\delta,\gamma}=\mathbb{C}\langle p_1,
p_2, \ldots, p_{k}, q_1, q_2,\ldots, q_{l}, s_1, s_2, \ldots, s_m |\\
 p_i p_j=\delta_{ij}p_i,
  q_i q_j= \delta_{ij}q_i; s_i s_j=
\delta_{ij}s_i;\\ \sum_{i=1}^{k}\alpha_i p_i
+\sum_{j=1}^{l}\beta_j q_j + \sum_{d=1}^{m}\delta_d s_d = \gamma e
  \rangle.
\end{gather*}
 Here $e$ is the identity of the algebra, $\delta_{ij}$ is the Kronecker symbol. This is a  $*$-algebra, if we declare all generators to be
self-adjoint. 

A $*$-representation $\pi$ of  $\mathcal{P}_{\alpha,\beta,\delta,\gamma}$ is determined by a  triple of non-negative operators $A=\sum_{i=1}^k \alpha_i P_i$, $B=\sum_{j=1}^l \beta_j Q_j$, $C=\sum_{d=1}^m \delta_d S_d$, where each of the  families of orthoprojections, $\{I-\sum_{i=1}^k P_i$, $P_i=\pi(p_i)$,  $i=1,\ldots, k  \}$,  $\{I-\sum_{j=1}^l Q_j$, and $Q_j=\pi(q_j)$, $j=1,\ldots, l \}$ and $\{I-\sum_{d=1}^m S_d$, $S_d=\pi(s_d)$, $d=1,\ldots, m \}$,   forms a resolution of the identity and  such that  $A+B+C=\gamma I$. 
So in terms of $*$-representations, the spectral  problem is a  problem consisting of the following two parts:   
1) a description of the set $\Sigma_{k,l,m}$ of the parameters $\alpha$, $\beta$, $\delta$, $\gamma$ for which there exist $*$-representations of $\mathcal{P}_{\alpha,\beta,\delta,\gamma}$;  
2) a description of *-representations $\pi$ of the *-algebra $ \mathcal{P}_{\alpha,\beta,\delta,\gamma}$.  
 A natural way to try to  solve the spectral  problem is to describe all irreducible $*$-representations up to unitary equivalence and then all *-representations as sums or direct integrals of irreducible representations. 

An essential step in this direction is to describe  of irreducible non-degenerate representations. Let us introduce non-degenerate representations and then reformulate the  strict spectral  problem in terms of *-representations of algebras.    
 Let us call a  $*$-representation  $\pi$  of the algebra
 $\mathcal{P}_{\alpha,\beta,\delta,\gamma}$  {\it non-degenerate},
if  $\pi(p_i)\not=0$, $\pi(q_j)\not=0$, $\pi(s_d)\not=0$  for
$1\le i\le
 k$, $1\le j\le l$, $1\le d\le m$, and  $\sum_{i=1}^{k}\pi(p_i)\not=I$,
 $\sum_{j=1}^{l}\pi(q_j)\not=I$, $\sum_{d=1}^{m}\pi(s_d)\not=I$.
  Consider the following sets: 
  $T_{k,l,m}=\{ ({\alpha},{\beta},  {\delta},\gamma) |  {\alpha}\in \mathbb{R}^{+k},{\beta}\in \mathbb{R}^{+l},  {\delta}\in \mathbb{R}^{+m}, \gamma\in  \mathbb{R^+}\}$ and   $W=\{ ({\alpha},{\beta},  {\delta},\gamma)\in T_{k,l,m} |$ { \it there is a non-degenerate *-representation of }  $\mathcal{P}_{\alpha,\beta,\delta,\gamma}\}$; they  depend only on   $(k,l,m)$. With an  integer vector    $(k,l,m)$ we will  associate non-oriented  star-shape  graph  $G$  with three branches of the lengths   $k$, $l$ and   $m$ stemming from  single root. Henceforth we will denote $W$ by  $W(G)$, and  $T_{k,l,m}$ by  $T(G)$, where $G$ is the tree mentioned above. We will need some more notations     $W_{irr}=\{ ({\alpha},{\beta},   {\delta},\gamma)\in T(G) |$    { \it there exists a non-degenerate irreducible *-representation of}    $\mathcal{P}_{\alpha,\beta,\delta,\gamma}\}$.
The strict spectral problem   for operators on a Hilbert space can be
reformulated in the following way: 1) for a given graph  $G$ describe the set  $W(G)$; 2) describe non-degenerate representations $\mathcal{P}_{\alpha,\beta,\delta,\gamma}$ up to unitary equivalence.  

\noindent 4. A common part of the spectral and strict spectral  problems is a  description of the set $W_{irr}(G)$ and  irreducible representations. The present article is  devoted to this common part for Dynkin graphs.   

If the  graph is a Dynkin graph, the problem is greatly simplified. The 
algebras  $\mathcal{P}_{\alpha,\beta,\delta,\gamma}$ associated
with Dynkin graphs  (res. extended Dynkin graphs) have a more simple
structure then in other cases. In particular, the  algebras
$\mathcal{P}_{\alpha,\beta,\delta,\gamma}$ are finite dimensional,
(res. have polynomial  growth) if and only if the associated graph is a
Dynkin graph (res. an extended Dynkin graph) (see.~\cite{melit}). As shown
in~\cite{roiter} irreducible representations of the algebras
associated with Dynkin graphs exist only in certain dimensions that are
bounded from above. In the paper we give a complete description of  $W(G)_{irr}$ for all Dynkin graphs $G$ and an algorithm for finding  all irreducible representations.

 In Sec.~\ref{s1} and  Sec.~2 we will show that the theory of $*$-representations of the algebras associated with such graphs can be reformulated in terms of locally-scalar graph representations.

   In Sec.~\ref{s2} we summarize the  calculation of possible parameters and   generalized dimensions of non-degenerate irreducible  $*$-representations of the  $*$-algebras associated with Dynkin graphs using the  machinery elaborated in~\cite{roiter}.

Since the description of   $*$-representations of the  $*$-algebras associated with Dynkin graphs   can be reduced to a  description of non-degenerate irreducible representations  (see~\cite{Zavod}),  in this article we essentially solve the  spectral  problem for the algebras associated with Dynkin graphs. 

The results of the paper were partially announced in~\cite{KPS}.

Let us remark that representations of quivers  were  also studied in connection with Horn's problem in~\cite{crawley}.

\section{
Locally-scalar graph representations and representations of the 
algebras generated by orthoprojections.}\label{s1}

Henceforth we will use definitions, notations and results
 about representations of graphs in category of
Hilbert spaces found in ~\cite{roiter}. Let us recall some of them.

A graph $G$ consists of a set of vertices $G_v$ a set of edges 
$G_e$ and a map $\varepsilon$ from $G_e$ into the set of
one- and two-element subsets of $G_v$ (the edge is mapped into the set of
incident vertices). Henceforth  we consider connected finite
graphs without cycles (trees). Fix a decomposition of $G_v$
 of the form   $G_v= {\overset{\circ}{G}}_v \sqcup
{\overset{\bullet}{G}}_v$ (unique up to permutation) such that for
each  $\alpha\in G_e$ one of the vertices from
$\varepsilon({\alpha})$ belongs to ${\overset{\circ}{G}}_v$ and  the
other to ${\overset{\bullet}{G}}_v$.  Vertices in
${\overset{\circ}{G}}_v$ will be called even, and those in the set
${\overset{\bullet}{G}}_v$ odd.
 Let us recall the definition of a  representation  $\Pi$ of a  graph $G$
 in the category of Hilbert spaces  $\mathcal{H}$. Let us associate with each vertex  $g\in G_v$  a Hilbert space
   $\Pi(g)= H_g\in \text{Ob}
\mathcal{H}$, and with each edge  $\gamma\in G_e$ such that
 $\varepsilon(\gamma) =\{g_1,g_2\}$  a pair of mutually adjoint
operators
  $\Pi(\gamma)=\{ \Gamma_{g_1,g_2}$, $\Gamma_{g_2,g_1}\}$, where
  $\Gamma_{g_1,g_2}:H_{g_2}\to H_{g_1}$. Construct a category
  $\rep (G,\mathcal{H})$, its objects are representations of the
graph  $G$  in  $\mathcal{H}$. A  morphism  $C: \Pi\to
\widetilde{\Pi}$ is a family   $\{C_g\}_{g\in G_v}$ of operators
  $C_g: \Pi(g)\to \widetilde{\Pi}(g)$ such that the following diagrams
  commute  for all edges   $\gamma_{g_2,g_1}\in G_e$:
\begin{equation*}
\begin{CD}
H_{g_1} @>\Gamma_{g_2,g_1}>>H_{g_2} \\  @V C_{g_1} VV  @VV
C_{g_2}V
\\ \widetilde{H}_{g_1}
@>\widetilde{\Gamma}_{g_2,g_1}>>\widetilde{H}_{g_2}
\end{CD}
\end{equation*}

 Let $M_g$ be the set of vertices connected with  $g$ by an edge.
 Let us define the  operators
\[
A_g= \sum_{g'\in M_g} \Gamma_{gg'}\Gamma_{g'g}.
\]

A representation    $\Pi$ in  $\rep (G,\mathcal{H})$  will be called
{\it locally-scalar}, if all operators $A_g$ are scalar, $A_g=
\alpha_g I_{H_g}$. The full subcategory $\rep (G,\mathcal{H})$,
  objects  of which are locally-scalar representations, will
be denoted by $\rep G$ and  called the category of
locally-scalar  representations of the   graph $G$.

Let us denote by  $V_{G}$ the real vector space consisting of sets
$x=(x_g)$ of real numbers $x_g$, $g\in G_v$. Elements $x$ of $V_G$
we will call  $G$-vectors. A vector $x=(x_g)$ is called positive, 
$x>0$, if $x\not=0$ and  $x_g\ge 0$ for all $g\in G_v$. Denote
$V_G^+= \{x\in V_G| x>0\}$.  If $\Pi$ is a finite dimensional
representation of the graph $G$ then the  $G$-vector $(d(g))$, where $d(g)=
\dim \Pi(g)$ is called the {\it dimension} of $\Pi$. If $A_g=f(g)
I_{H_g}$ then the $G$-vector $f=(f(g))$ is called a {\it character } of a
locally-scalar representation   $\Pi$ and $\Pi$ is called an 
$f$-representation in this case. The {\it support} $G_v^{\Pi}$ of
 $\Pi$ is $\{ g\in G_v| \Pi(g) \not= 0 \}$. A representation   $\Pi$
is {\it faithful} if $G_v^{\Pi}= G_v$. A character of the
locally-scalar representation $\Pi$ is uniquely defined on the
support  $G_v^{\Pi}$ and non-uniquely on its complement. In the 
general case, denote by $\{f_\Pi \}$ the set of characters of
$\Pi$. For each vertex $g\in G_v$, denote by  $\sigma_g$ the linear
operator on  $V_G$ given by the formulae:
 \begin{gather*}
 (\sigma_g x)_{g'} = x_{g'}\ \text{if}\ g'\not=g,\\
(\sigma_g x)_{g} = -x_{g} +\sum_{g'\in M_g}x_{g'}.
 \end{gather*}
The mapping   $\sigma_g$  is called a {\it reflection} at the vertex
$g$. The composition of all reflections at odd vertices is denoted
by   $ \overset{\bullet}{c}$ (it does not depend on the order of the 
factors), and at all even vertices by $ \overset{\circ}{c}$.
A Coxeter transformation is $c=\overset{\circ}{c}\overset{\bullet}{c}$,
$c^{-1}= \overset{\bullet}{c}\overset{\circ}{c}$.  The transformation
$\overset{\bullet}{c}$ ($\overset{\circ}{c}$)  is called an  odd
(even) Coxeter map.
  Let us adopt the  following notations for compositions of the  Coxeter maps: 
   $\overset{\bullet}{c}_k= \ldots
\overset{\bullet}{c}\overset{\circ}{c}\overset{\bullet}{c}$ ($k$-
factors), $\overset{\circ}{c}_k= \ldots
\overset{\circ}{c}\overset{\bullet}{c}\overset{\circ}{c}$ ($k$-
factors), $k\in\mathbb{N}$.

Any real function  $f$ on $G_v$ can be identified  with a $G$-vector
$f=(f(g))_{g\in G_v}$. If $d(g)$  is the dimension of a
locally-scalar graph representation  $\Pi$, then
\begin{gather}
\overset{\circ}{c}(d)(g)=\begin{cases} -d(g) +  \sum_{g'\in
M_g}d(g'), &\text{if}\  g\in \overset{\circ}{G}_v,\\
d(g), & \text{if}\ g\in \overset{\bullet}{G}_v,
\end{cases}\\
\overset{\bullet}{c}(d)(g)=\begin{cases} -d(g) +  \sum_{g'\in
M_g}d(g'), &\text{if}\  g\in \overset{\bullet}{G}_v,\\
d(g), & \text{if}\ g\in \overset{\circ}{G}_v.
\end{cases}
\end{gather}

 For $ d\in Z_G^+$ and  $f\in V_G^+$, consider the full subcategory
  $\rep (G,d,f)$ in $\rep
 G$ (here $Z_G^+$ is the set of positive integer  $G$-vectors),
with the set of objects  $Ob \rep (G, d, f)= \{ \Pi| \dim \Pi (g)
=d(g), f\in\{f_\Pi\}\}$. All representations
    $\Pi$ from $\rep (G,d,f)$ have the same support  $X=X_d=G_v^\Pi=\{g\in G_v|
  d(g)\not=0\}$.   We will consider these categories only if   $(d,f)\in
  S=\{ (d,f)\in Z_G^+\times V_G^+ | d(g) + f(g) >0, g\in G_v \}$.
   Let  $\overset{\circ}{X}= X\cap \overset{\circ}{G}_v$,
    $\overset{\bullet}{X}= X\cap \overset{\bullet}{G}_v$.  $\rep_{\circ}(G, d, f) \subset \rep(G, d,
    f)$  ( $\rep_{\bullet}(G, d, f) \subset \rep(G, d,
    f)$) is the full subcategory with objects  $(\Pi, f)$ where
      $f(g)>0$  if  $g\in \overset{\circ}{X}$ ($f(g)>0$  if  $g\in
    \overset{\bullet}{X}$). Let $S_0=\{(d,f)\in S|f(g)>0 \text{ if }
g\in \overset{\circ}{X}_d\}$, $S_\bullet=\{(d,f)\in S|f(g)>0
\text{ if } g\in \overset{\bullet}{X}_d\}$

Put
\begin{gather}
\overset{\bullet}{c}_d(f)(g)=\overset{\circ}{f}_d(g)=
\begin{cases}
 \overset{\bullet}{c}(f)(g), &\text{if}\  g\in \overset{\bullet}{X}_d,\\
f(g), & \text{if}\ g\not\in \overset{\bullet}{X}_d,
\end{cases}\\
\overset{\circ}{c}_d(f)(g)=\overset{\bullet}{f}_d(g)=
\begin{cases}
 \overset{\circ}{c}(f)(g), &\text{if}\  g\in \overset{\circ}{X}_d,\\
f(g), & \text{if}\ g\not\in \overset{\circ}{X}_d.
\end{cases}.
\end{gather}

 Let us denote  $\overset{\bullet}{c}_d^{(k)}(f)=
 \ldots \overset{\bullet}{c}_{\overset{\circ}{c}_2(d)}
 \overset{\circ}{c}_{\overset{\circ}{c}(d)}
 \overset{\bullet}{c}_d(f)$ ($k$  factors), 
 $\overset{\circ}{c}_d^{(k)}(f)=
 \ldots \overset{\circ}{c}_{\overset{\bullet}{c}_2(d)}
  \overset{\bullet}{c}_{\overset{\bullet}{c}(d)}
  \overset{\circ}{c}_d(f)$ ($k$  factors).
The even and odd Coxeter reflection functors are defined
in~\cite{roiter}, 
  $\overset{\circ}{F}: \rep_{\circ} (G, d, f) \to
  \rep_{\circ} (G, \overset{\circ}{c}(d), \overset{\circ}{f}_d)$
   if  $(d,f)\in S_\circ$,
$\overset{\bullet}{F}: \rep_{\bullet} (G, d, f) \to
  \rep_{\bullet} (G, \overset{\bullet}{c}(d), \overset{\bullet}{f}_d)$
   if  $(d,f)\in S_\bullet$; they are equivalences of the  categories.
Let us denote
   $\overset{\circ}{F}_k(\Pi)=
 \ldots \overset{\circ}{F}\overset{\bullet}{F}\overset{\circ}{F}(\Pi)$ ($k$  factors), 
 $\overset{\bullet}{F}_k(\Pi)=
 \ldots \overset{\bullet}{F}\overset{\circ}{F}\overset{\bullet}{F}(\Pi)$ ($k$
 factors),  if the compositions exist. Using these functors, an analog of Gabriel's theorem
for graphs and their locally-scalar representations has been
proven in~\cite{roiter}. In particular, it has been  proved that
any locally-scalar graph representation decomposes into a direct sum
(finite or infinite) of finite dimensional indecomposable
representations, and all indecomposable representations can be
obtained by odd and even Coxeter reflection functors starting from
the simplest representations  $\Pi_g$ of the graph  $G$
 ($\Pi_g(g)=\mathbb{C}, \Pi_g(g')= 0$  if $g\not=g';\ g,g'\in G_v
 $).

\section{ Representations of algebras generated by projections.}

Let us consider a tree $G$ with vertices  $\{ g_i,i=0,\ldots,
k+l+m\}$ and edges $\gamma_{g_i g_j}$,  see the figure.

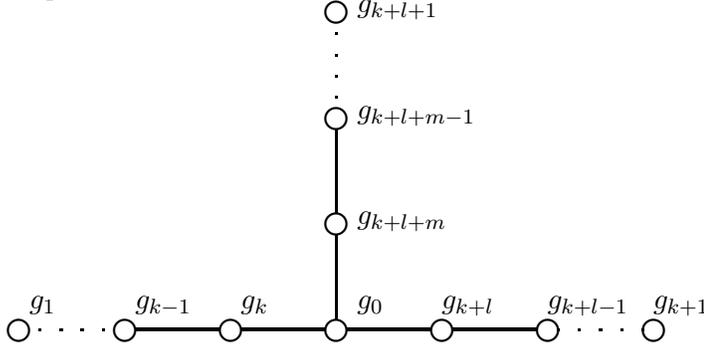
\setlength{\unitlength}{4pt}
\begin{picture}(70,32)(-1,-1)
\thicklines
 \drawline(29,1)(29,9)
\drawline(29,11)(29,19)
 \multiputlist(0,0)(10,0){\circle{2},\circle{2},\circle{2},\circle{2},\circle{2},\circle{2},\circle{2}}
 \multiputlist(30,10)(0,10){\circle{2},\circle{2},\circle{2}}
 \drawline(10,0)(18,0)
\drawline(20,0)(28,0) \drawline(30,0)(38,0)\drawline(40,0)(48,0)
\dottedline{2}(1,0)(7,0) \dottedline{2}(50,0)(58,0)
\dottedline{2}(29,22)(29,28)
 \put(31,10){$g_{k+l+m}$}\put(31,20){$g_{k+l+m-1}$}\put(31,30){$g_{k+l+1}$}
 \put(39,2){$g_{k+l}$}\put(49,2){$g_{k+l-1}$}\put(59,2){$g_{k+1}$}
\put(0,2){$g_{1}$}\put(10,2){$g_{k-1}$}\put(20,2){$g_{k}$}\put(31,2){$g_{0}$}
\end{picture}
\vspace{10pt}

We will establish a connection  between  $*$-representations of the
$*$-algebra  $\mathcal{P}_{\alpha,\beta,\delta,\gamma}$  and
locally-scalar representations of the graph  $G$, see~\cite{KPS},  $\alpha=(\alpha_1,\alpha_2,\ldots, \alpha_k)$, $\beta =(\beta_1, \beta_2,\ldots, \beta_l)$, $\delta=(\delta_1,\delta_2,\ldots, \delta_m)$.

\begin{definition}
An irreducible finite dimensional   *-representation  $\pi$ of the
algebra  $
 \mathcal{P}_{\alpha,\beta,\delta,\gamma}$ such that
 $\pi(p_i)\not= 0 (1\le i\le k),\ \pi(q_j)
 \not= 0 (1\le j\le l),
 \ \pi(s_d)\not= 0 (1\le d\le m)$, and
 $\sum_{i=1}^k \pi(p_i)\not= I, \sum_{j=1}^l \pi(q_j)\not= I,
 \sum_{d=1}^m \pi(s_d)\not= I$,  will be called
  {\it non-degenerate}. By
 $\overline{\rep} \mathcal{P}_{\alpha,\beta,\delta,\gamma}$
 we will denote the full subcategory of non-degenerate
representations in the category ${\rep}
\mathcal{P}_{\alpha,\beta,\delta,\gamma}$ of $*$-representations of the $*$-algebra $\mathcal{P}_{\alpha,\beta,\delta,\gamma}$ in the
category $\mathcal{H}$ of Hilbert spaces.
\end{definition}

Let $\pi$ be a $*$-representation of
$\mathcal{P}_{\alpha,\beta,\delta,\gamma}$
 on a Hilbert space  $H_0$. Put
  $P_i=\pi(p_i)$, $1\le i \le k$,
 $Q_j=\pi(q_j)$, $1\le j \le l$, $S_t=\pi(s_t)$, $1\le t \le m$.
 Let   $H_{p_i}= \image P_i$, $H_{q_j}= \image Q_j$, $H_{s_t}=
 \image S_t$. Denote by  $\Gamma_{p_i}$, $\Gamma_{q_j}$,
 $\Gamma_{s_t}$ the corresponding  natural isometries
 $H_{p_i}\to H_0$, $H_{q_j}\to H_0$,  $H_{s_t}\to H_0$.
 Then, in particular,  $\Gamma_{p_i}^*\Gamma_{p_i}
 =I_{H_{p_i}}$ is the  identity operator on  $H_{p_i}$ and   $\Gamma_{p_i}\Gamma_{p_i}^*
 =P_i$. Similar equalities hold for the operators
 $\Gamma_{q_i}$ and $\Gamma_{s_i}$. Using  $\pi$  we construct a locally-scalar
representation   $\Pi$ of the graph   $G$.

  Let   $\Gamma_{ij}:H_j\to H_i$ denote the operator adjoint to   $\Gamma_{ji}:H_i\to H_j$, i.e. $\Gamma_{ij}=\Gamma_{ji}^*$. Put 
\begin{gather*}
\Pi(g_0)= H^{g_0}=H_0,\\
 \Pi(g_k)= H^{g_k}= H_{p_1}\oplus H_{p_2}\oplus \ldots \oplus
H_{p_k},\\
\Pi(g_{k-1})= H^{g_{k-1}}=  H_{p_2}\oplus \ldots \oplus
H_{p_{k-1}}\oplus H_{p_k},\\
\Pi(g_{k-2})= H^{g_{k-2}}= H_{p_2}\oplus H_{p_3}\oplus \ldots
\oplus
H_{p_{k-1}},\\
\ldots.
 \end{gather*}

In these equalities the summands are omitted  from the left and the right in turns. Analogously, we define subspaces  $\Pi(g_i)$ for $i=k+1,\ldots,k+l$ and  $i=k+l+1,\ldots,k+l+m$. Define  the
 operators  $\Gamma_{g_0,g_i}:H^{g_i}\to H^{g_0}$, where
$i\in\{k,k+l,k+l+m\}$, by the block-diagonal matrices
\begin{gather*}
\Gamma_{g_0,g_k}=\big[ \sqrt{\alpha_1} \Gamma_{p_1}|
\sqrt{\alpha_2} \Gamma_{p_2}|\ldots | \sqrt{\alpha_k} \Gamma_{p_k}
\big],\\
\Gamma_{g_0,g_{k+l}}=\big[ \sqrt{\beta_1} \Gamma_{q_1}|
\sqrt{\beta_2} \Gamma_{q_2}|\ldots | \sqrt{\beta_l} \Gamma_{q_l}
\big],\\
\Gamma_{g_0,g_{k+l+m}}=\big[ \sqrt{\delta_1} \Gamma_{s_1}|
\sqrt{\delta_2} \Gamma_{s_2}|\ldots | \sqrt{\delta_m} \Gamma_{s_m}
\big].\\
\end{gather*}
Now we define the representation  $\Pi$ on the edges  $\gamma_{g_0,g_k}$, $\gamma_{g_0,g_{k+l}}$, $\gamma_{g_0,g_{k+l+m}}$ by the rule
\begin{gather*}
\Pi(\gamma_{g_0,g_k})=\{\Gamma_{g_0,g_k},\Gamma_{g_k,g_0} \},\\
\Pi(\gamma_{g_0,g_{k+l}})=\{\Gamma_{g_0,g_{k+l}},\Gamma_{g_{k+l},g_0} \},\\
\Pi(\gamma_{g_0,g_{k+l+m}})=\{\Gamma_{g_0,g_{k+l+m}},\Gamma_{g_{k+l+m},g_0} \}.\\
\end{gather*}

It is easy to see that 
\[\Gamma_{g_0,g_k}
\Gamma_{g_k,g_0}+\Gamma_{g_0,g_{k+l}}
\Gamma_{g_{k+l},g_0}+\Gamma_{g_0,g_{k+l+m}}
\Gamma_{g_{k+l+m},g_0}=\gamma I_{H^{g_0}}.
\]

Let  $\mathcal{O}_{H,0}$  denote the operators from the zero space to  $H$, and  $\mathcal{O}_{0,H}$  denote the  zero operator from  
$H$ into the zero subspace. For the  operators $\Gamma_{g_j,g_i}
:H^{g_i}\to H^{g_j}$ with $i,j\not=0$, put
\begin{multline}\label{zved}
\Gamma_{g_{k-1},g_k}=\mathcal{O}_{0,H_{p_1}}\oplus
\sqrt{\alpha_{1} - \alpha_2} I_{H_{p_2}} \oplus \sqrt{\alpha_{1} -
\alpha_3} I_{H_{p_3}}\oplus \ldots \oplus \sqrt{\alpha_{1} -
\alpha_{k}}
I_{H_{p_{k}}}, \\
 \Gamma_{g_{k-1},g_{k-2}}=\sqrt{\alpha_{2}
- \alpha_k} I_{H_{p_2}} \oplus \sqrt{\alpha_{3} - \alpha_k}
I_{H_{p_3}} \oplus \ldots \oplus \sqrt{\alpha_{k-1} - \alpha_{k}}
I_{H_{p_{k-1}}}\oplus\mathcal{O}_{H_{p_k},0},\\
\Gamma_{g_{k-3},g_{k-2}}=\mathcal{O}_{0,H_{p_{2}}}\oplus
\sqrt{\alpha_{2} - \alpha_3} I_{H_{p_3}} \oplus \sqrt{\alpha_{2} -
\alpha_4} I_{H_{p_4}}\oplus \ldots \oplus \sqrt{\alpha_{2} -
\alpha_{k-1}}
I_{H_{p_{k-1}}},\\
\ldots\ldots\ldots\ldots\ldots\ldots\ldots\ldots\ldots\ldots.
\end{multline}

The corresponding operators for the rest of the edges of 
 $G$ can be constructed analogously. One can
check that the  operators  $\Gamma_{g_{i},g_{j}}$, where
$\Gamma_{g_{i},g_{j}}=\Gamma_{g_{i},g_{j}}^*$, define a  locally-scalar representation of the graph  $G$ with the following
character  $f$:

\begin{gather}\label{char}
\begin{alignat*}{3}
 f(g_k) &= \alpha_1, &   \qquad f(g_{k+l}) &=
\beta_1, &  \qquad f(g_{k+l+m}) &= \delta_1, \\
 f(g_{k-1}) &=
\alpha_{1}-\alpha_{k}, & \qquad f(g_{k+l-1}) &=
\beta_{1}-\beta_{l}, &
\qquad f(g_{k+l+m-1})&= \delta_{1}-\delta_{m}, \\
 f(g_{k-2}) &= \alpha_{2}-\alpha_{k}, & \qquad f(g_{k+l-2}) &=
\beta_{2}-\beta_{l}, & \qquad f(g_{k+l+m-2}) &=
\delta_{2}-\delta_{m}, \\
   f(g_{k-3}) &= \alpha_{2}-\alpha_{k-1}, & \qquad f(g_{k+l-3}) &=
\beta_{2}-\beta_{l-1}, & \qquad f(g_{k+l+m-3}) &=
\delta_{2}-\delta_{m-1}, \\
 f(g_{k-4}) &= \alpha_{3}-\alpha_{k-1}, & \qquad f(g_{k+l-4}) &=
\beta_{3}-\beta_{l-1}, & \qquad f(g_{k+l+m-4}) &=
\delta_{3}-\delta_{m-1},\\
  \ldots &   & \ldots  & &  \ldots \\
\end{alignat*}
\end{gather}

$f(g_0)=\gamma$. 

And vice versa, if a locally-scalar representation of the graph $G$ 
with the  character $f(g_i)=x_i\in\mathbb{R}^*$ corresponds to a non-degenerate representation of  $\mathcal{P}_{\alpha,\beta,\delta,\gamma}$, then one can check that  
\begin{gather*}
\alpha_1=x_k,\\
\alpha_k=x_k-x_{k-1},
\alpha_{2}=x_k-x_{k-1}+x_{k-2},\\
\alpha_{k-1}=x_k-x_{k-1}+x_{k-2}-x_{k-3},\\
\alpha_{3}=x_k-x_{k-1}+x_{k-2}-x_{k-3}+x_{k-4},\\
\ldots.
\end{gather*}
Here $x_j=0$ if $j\le 0$. 
Analogously one can find   $\beta_j$ and  $\delta_t$. We will denote   $\Pi$ by $\Phi(\pi)$.

 Let  $\pi$  and $\widetilde{\pi}$ be non-degenerate representations of the algebra $\mathcal{P}_{\alpha,\beta,\delta,\gamma}$ and $C_0$ an intertwining operator for these representations; this is a morphism from  $\pi$ to  $\widetilde{\pi}$ in the category  $\rep G$), 
 $C_0: H_0\to \widetilde{H}_0$,  $C_0 \pi= \widetilde{\pi} C_0$.  Put
\begin{gather*}
C_{p_i}= \widetilde{\Gamma}_{p_i}^* C_0 \Gamma_{p_i}, C_{p_i}:
H_{p_i}\to \widetilde{H}_{p_i},\   1\le i \le k,\\
C_{q_j}= \widetilde{\Gamma}_{q_j}^* C_0 \Gamma_{q_j}, C_{q_j}:
H_{q_j}\to \widetilde{H}_{q_j},\   k+1\le j \le k+l,\\
C_{s_t}= \widetilde{\Gamma}_{s_t}^* C_0 \Gamma_{s_t}, C_{s_t}:
H_{s_t}\to \widetilde{H}_{s_t}.\   k+l+1\le t \le k+l+m,\\
\ldots
\end{gather*}
 Put
\begin{gather*}
C^{(g_0)}=C_0: H^{(g_0)}\to \widetilde{H}^{(g_0)},\\
 C^{(g_k)}= C_{p_1}\oplus C_{p_2} \oplus \ldots  \oplus C_{p_k} :
H^{(g_k)}\to \widetilde{H}^{(g_k)},\\
 C^{(g_{k-1})}=  C_{p_2} \oplus \ldots  \oplus C_{p_{k-1}} \oplus C_{p_k}:
H^{(g_{k-1})}\to \widetilde{H}^{(g_{k-1})},\\
 C^{(g_{k-2})}= C_{p_2}\oplus C_{p_3} \oplus \ldots  \oplus C_{p_{k-1}} :
H^{(g_{k-2})}\to \widetilde{H}^{(g_{k-2})},\\
\ldots
\end{gather*}
 Analogously one can construct the operators  $C^{(g_i)}$ for  $i\in \{ k+l,\ldots, k+l+m
 \}$.  It is routine to check that the operators  $\{C^{(g_i)}\}_{0\le i\le
 k+l+m}$ intertwine the representations   $\Pi=\Phi(\pi)$  and 
 $\widetilde{\Pi}=\Phi(\widetilde{\pi})$.  Put    $\Phi(C_0)=\{ C^{(g_i)}_{0\le i\le
 k+l+m}\}$.  Thus we have defined a functor   $\Phi:
 \overline{\rep} \mathcal{P}_{\alpha,\beta,\delta,\gamma}\to
 {\rep} G$, see~\cite{KPS}.   Moreover, the functor   $\Phi$  is univalent and full.
Let  $\widetilde{\rep}(G,d,f)$ be the full subcategory of irreducible representations of  ${\rep}(G,d,f)$.
   $\Pi\in \ob \widetilde{\rep}(G,d,f)$, $f(g_i)= x_i\in \mathbb{R}^+,
 d(g_i)=d_i\in \mathbb{N}_0$, where $f$ is  the character of   $\Pi$, $d$ its dimension. It easy to verify that the representation $\Pi$ is isomorphic  (unitary equivalent) to an irreducible representation from the image of the functor $\Phi$ if and only if 
\begin{gather}\label{nondeg}
1.\    0<x_1<x_2<\ldots < x_k; 0<x_{k+1}<x_{k+2}<\ldots <
x_{k+l};\\
  0<x_{k+l+1}<x_{k+l+2}<\ldots < x_{k+l+m}; \\
2.\    0< d_1<d_2<\ldots < d_k< d_0;  0< d_{k+1}<d_{k+2}<\ldots < \label{ineq1}\\
d_{k+l}< d_0;  0< d_{k+l+1}<d_{k+l+2}<\ldots < d_{k+l+m}<
d_0.\label{ineq2}
\end{gather}
(All matrices of the representation of the graph $G$,  except for 
$\Gamma_{g_0,g_k}, \Gamma_{g_k,g_0}, \Gamma_{g_0,g_{k+l}}$, $
\Gamma_{g_{k+l}, g_0}$, $\Gamma_{g_0,g_{k+l+m}},
\Gamma_{g_{k+l+m}, g_0}$, can be brought    to the "canonical" form (\ref{zved}) by admissible transformations. Then the rest of the matrices will naturally be  partitioned into  blocks, which gives the matrices $\Gamma_{p_i}, \Gamma_{q_i}, \Gamma_{s_i}$, and thus the projections $P_i, Q_i, R_i$).  An irreducible representation  $\Pi$ of the graph $G$ satisfying conditions  (\ref{nondeg})--(\ref{ineq2})  will be called {\it non-degenerate}.
 Let 
\begin{gather*}
\dim H_{p_i}= n_i, 1\le i\le k; \\
\dim H_{q_j}= n_{k+j}, 1\le j\le l; \\
\dim H_{s_t}= n_{k+l+t}, 1\le t\le m; \\
\dim H_{0}= n_{0}.
\end{gather*}
The vector ${n}= (n_0,n_1,\ldots , n_{k+l+m})$ is called a {\it  generalized dimension} of the representation  $\pi$ of the algebra 
$\mathcal{P}_{\alpha,\beta,\delta,\gamma}$.  Let $\Pi=\Phi(\pi)$
for a  non-degenerate representation of the algebra 
$\mathcal{P}_{\alpha,\beta,\delta,\gamma}$, ${d}=(d_0,d_1,\ldots,
d_{k+l+m})$ be the dimension of  $\Pi$.  It is easy to see that 
\begin{gather*}
n_1+n_2+\ldots + n_k = d_k,\\
n_2+\ldots + n_{k-1}+n_k = d_{k-1},\\
    n_2+\ldots + n_{k-1} = d_{k-2},\\
    n_3+\ldots + n_{k-2}+n_{k-1} = d_{k-3},\\
    \ldots
\end{gather*}
 Thus 
\begin{multline}\label{dim}\\
n_1 =d_k-d_{k-1},\\
n_k =d_{k-1}-d_{k-2},\\
n_{2} =d_{k-2}-d_{k-3},\\
    \ldots \\
\end{multline}
  Analogously one can find  $n_{k+1}, \ldots , n_{k+l}$ from 
  $d_{k+1}, \ldots, d_{k+l}$  and  $n_{k+l+1}, \ldots ,n_{k+l+m}$ from  $d_{k+l+1}, \ldots, d_{k+l+m}$

 Denote by  $\overline{\rep}G$  the full subcategory in   $\rep
 G$ of non-degenerate locally-scalar representations of the graph   $G$.
As a corollary of the previous arguments we obtain the following theorem. 
\begin{theorem}
Let  $\mathcal{P}_{\alpha,\beta,\delta,\gamma}$ be associated with a
graph  $G$. The functor  $\Phi$ is an equivalence of the categories
$\overline{\rep}{ \mathcal{P}_{\alpha,\beta,\delta,\gamma}}$
 of  non-degenerate *-representations of the algebra
 $\mathcal{P}_{\alpha,\beta,\delta,\gamma}$ and the
category $\overline{\rep}G$ of non-degenerate locally-scalar
representations of the graph  $G$.
\end{theorem}

Let us define the Coxeter functors for the  $*$-algebras 
$\mathcal{P}_{\alpha,\beta,\delta,\gamma}$, by putting 
$\overset{\circ}{\Psi}= \Phi^{-1}\overset{\circ}{F} \Phi $ and 
$\overset{\bullet}{\Psi}= \Phi^{-1}\overset{\bullet}{F} \Phi$. 
Now we can use the results of~\cite{roiter} to give a description of representation of the *-algebra 
$\mathcal{P}_{\alpha,\beta,\delta,\gamma}$.

Obviously, all irreducible representations of the *-algebra  associated with the Dynkin diagram $A_n$ are one-dimensional, with diagram  $D_n$ one- or two-dimensional  and only the  algebras associated with the  diagrams $D_4, E_6, E_7, E_8$ have non-degenerate irreducible representations.

In the next section we will do the following: we know how to construct all irreducible locally-scalar representations of Dynkin graphs with the aid of Coxeter reflection functors starting from the simplest ones. In particular, we can find their dimensions and characters \cite{roiter}. Next we single out non-generate representations and apply the  equivalence functor  $\Phi$, see  Theorem~1.

\section{Algebras associated with Dynkin graphs.}\label{s2}
Let  $G$ be the Dynkin graph $E_6$, 

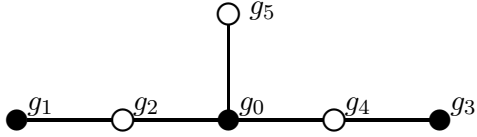
\begin{picture}(62,15)(-1,-1)
\thicklines
 \drawline(19,0)(19,9)
 \multiputlist(0,0)(10,0){\circle*{2},\circle{2},\circle*{2},\circle{2},\circle*{2}}
 \multiputlist(20,10)(0,10){\circle{2}}
 \drawline(0,0)(8,0)
\drawline(10,0)(19,0) \drawline(20,0)(28,0)\drawline(30,0)(38,0)
\put(0,1){$g_1$} \put(10,1){$g_2$} \put(20,1){$g_0$}
\put(30,1){$g_4$} \put(40,1){$g_3$} \put(21,10){$g_5$}
\end{picture}

 The vertices   $g_0, g_1, g_3$ will be called odd and marked with    $\bullet$ on the graph, the vertices $g_2, g_4, g_5$  are even and indicated with $\circ$. The parameters of the  corresponding algebra  $\mathcal{P}_{\alpha,\beta,\delta,\gamma}$ are enumerated according to the following picture:
$$
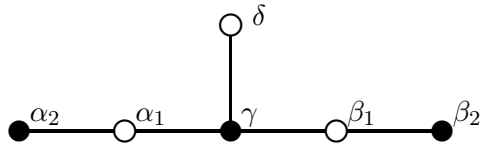
\begin{picture}(62,15)(-1,-1)
\thicklines
 \drawline(19,0)(19,9)
 \multiputlist(0,0)(10,0){\circle*{2},\circle{2},\circle*{2},\circle{2},\circle*{2}}
 \multiputlist(20,10)(0,10){\circle{2}}
 \drawline(0,0)(8,0)
\drawline(10,0)(19,0) \drawline(20,0)(28,0)\drawline(30,0)(38,0)
\put(0,1){$\alpha_2$} \put(10,1){$\alpha_1$} \put(20,1){$\gamma$}
\put(30,1){$\beta_1$} \put(40,1){$\beta_2$} \put(21,10){$\delta$}
\end{picture}
$$

Using  the proof of  theorem~3.9 (see \cite{roiter}) and  direct calculations we see that all non-degenerate representations of the graph  $E_6$ can be obtained by the Coxeter reflection  functors   $\overset{\circ}{F}$ and  $\overset{\bullet}{F}$ in   4 and  5 steps  from the simplest one $\Pi_{g_0}$ that corresponds to the simple root   $\overline{g}_0$  of the graph   $E_6$. The dimension of the representation   $\Pi_{g_0}$ is denoted by $d_{g_0}(g)$.  Then $d_{g_0}(g_i)=0$ if $g_i\not=g_0$,  $d_{g_0}(g_0)=1$. Let the character of   $\Pi_{g_0}$ is denoted by  $f_{g_0}$,  then  $f_{g_0}(g_0)=0$, $f_{g_0}(g_i)=t_i>0$  if $g_i\not=g_0$.  All irreducible locally-scalar representations of Dynkin diagrams are obtained from the simple representations with characters of such type  by using the Coxeter reflection functors  (see~\cite{roiter}).  In the subsequent  calculations we will write the character value below or to the right  and the dimension  over or to the left of the vertex.    
$$
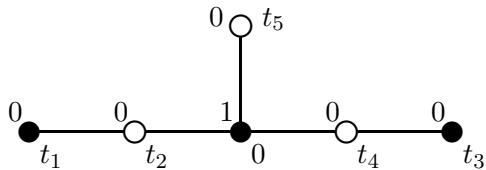
\begin{picture}(62,20)(-1,-5)
\thicklines
 \drawline(19,0)(19,9)
 \multiputlist(0,0)(10,0){\circle*{2},\circle{2},\circle*{2},\circle{2},\circle*{2}}
 \multiputlist(20,10)(0,10){\circle{2}}
 \drawline(0,0)(8,0)
\drawline(10,0)(19,0) \drawline(20,0)(28,0)\drawline(30,0)(38,0)
\put(-3,1){$0$} \put(7,1){$0$} \put(17,1){$1$} \put(27,1){$0$}
\put(37,1){$0$} \put(16,10){$0$} \put(0,-3){$t_1$}
\put(10,-3){$t_2$} \put(20,-3){$0$} \put(30,-3){$t_4$}
\put(40,-3){$t_3$} \put(21,10){$t_5$}
\end{picture}
$$

 Thus  $ \Pi= \overset{\circ}{F}_4 (\Pi_{g_0})$ is a non-degenerate representation of the graph   $G$ with dimension  $d= \overset{\circ}{c}_4(d_{g_0})$.  Calculations  values of the dimensions are summarized in the following table: 
$$
\begin{tabular}{|c|c|c|c|c|c|c|}
\hline
$g_i$ & $g_0$ & $g_1$ & $g_2$ & $g_3$ & $g_4$ & $g_5$  \\
  \hline
 $d$&$3$ & $1$ & $2$ & $1$ & $2$ & $1$
  \\
 \hline
\end{tabular}
$$
Then $\Pi_{g_0}= \overset{\bullet}{F}_4(\Pi)$, $d_{g_0}=
\overset{\bullet}{c}_4(d)$ і
 $f_{g_0}(g)=\overset{\circ}{c}_d^{(4)}(f)$.

Thus \begin{equation}\label{case}
\begin{cases}
t_i= \overset{\circ}{c}_d^{(4)}(f)(g_i)>0,  &\text{if}\  g_i\not=g_0,\\
t_0= \overset{\circ}{c}_d^{(4)}(f)(g_0)=0,  &\text{if}\  g=g_0.
\end{cases}
\end{equation}

 Let $\pi$ be a non-degenerate representation of the algebra   
  $\mathcal{P}_{\alpha,\beta,\delta,\gamma}$.  Then  $\Phi(\pi)$
  is non-degenerate representation from   $\overline{\rep} G$  with the  character   $f$ given by the following table  (see equalities~(\ref{char})):
$$
\begin{tabular}{|c|c|c|c|c|c|c|}
\hline
$g_i$ & $g_0$ & $g_1$ & $g_2$ & $g_3$ & $g_4$ & $g_5$  \\
  \hline
$f$ & $\gamma$ &
 $\alpha_1-\alpha_2$ & $\alpha_1$ & $\beta_1-\beta_2$ & $\beta_1$ & $\delta$
  \\
 \hline
\end{tabular}
$$

Conditions~(\ref{case}) imposed on $f(g_i)$ are necessary and sufficient for a representation of the algebra  
$\mathcal{P}_{\alpha,\beta,\delta,\gamma}$ in the dimension 
$\Phi(d)$ to exist. Computing the character we obtain 
$\overset{\circ}{c}_d^{(4)}(f)$, and from~(\ref{case}) we get the following inequalities: 
\[
\gamma>\alpha_1, \gamma>\beta_1, \alpha_2+\beta_1>\gamma,\\
\alpha_1+\beta_2>\gamma, \gamma>\alpha_2+\beta_2
\]
and the identity    $3 \gamma =
\delta+\alpha_1+\alpha_2+\beta_1+\beta_2$. Knowing the dimension  $d$  we find the corresponding generalized dimension   
$\Phi^{-1}(d)= (n_{g_0}; n_{g_1},\ldots, n_{g_5}) $ of the representation $\pi$,  see formulae~(\ref{dim}), where   $n_{g_i}=n_i$.
$$
\begin{tabular}{|c|c|c|c|c|c|c|}
\hline
$g_i$ &$g_0$ & $g_1$ & $g_2$ & $g_3$ & $g_4$ & $g_5$  \\
  \hline
 $\Phi^{-1}(d)$ & 3 & 1 & 1 & 1 & 1 & 1
  \\
 \hline
\end{tabular}
$$

In a similar way,  one can perform calculations for the representation   
 $\overset{\circ}{F}_5(\Pi_{g_0})$,  the only remaining non-degenerate irreducible representation of the  Dynkin graph   $E_6$. Thus we get the following theorem.  
\begin{theorem}\label{ineqe6}
 The algebra  $\mathcal{P}_{(\alpha_1, \alpha_2 ),
(\beta_1, \beta_2),\delta,\gamma}$  associated with the  Dynkin graph  $E_6$ has an irreducible non-degenerate representation only in the following cases. 
\begin{enumerate}
\item The following conditions are satisfied: 
\[
\gamma>\alpha_1, \gamma>\beta_1, \alpha_2+\beta_1>\gamma,\\
\alpha_1+\beta_2>\gamma, \gamma>\alpha_2+\beta_2,
\] and $3 \gamma =
\delta+\alpha_1+\alpha_2+\beta_1+\beta_2$; in this case there is a unique  representation in the  generalized dimension   $(3;1,1,1,1,1)$
\item
The following conditions are satisfied:
\[
\delta+\beta_1>\gamma, \gamma>\delta+\alpha_2,
\gamma>\delta+\beta_2, \delta+\alpha_1>\gamma,
\delta+\alpha_2+\beta_2>\gamma,
\] and  $3 \gamma =
2 \delta+\alpha_1+\alpha_2+\beta_1+\beta_2$; in this case there is a unique  representation in the  generalized dimension   $(3;1,1,1,1,2)$
\end{enumerate}
\end{theorem}

We have shown in the example of  $E_6$ how to perform calculations and  derive conditions of existence for an irreducible non-degenerate representation of the corresponding $*$-algebra. 
Now without additional explanations we  present complete answers for the rest of the Dynkin graphs. We label vertices    $g_i$ of the graph  $G$  as in the following pictures and write the corresponding 
  parameters of the algebra in parentheses: 

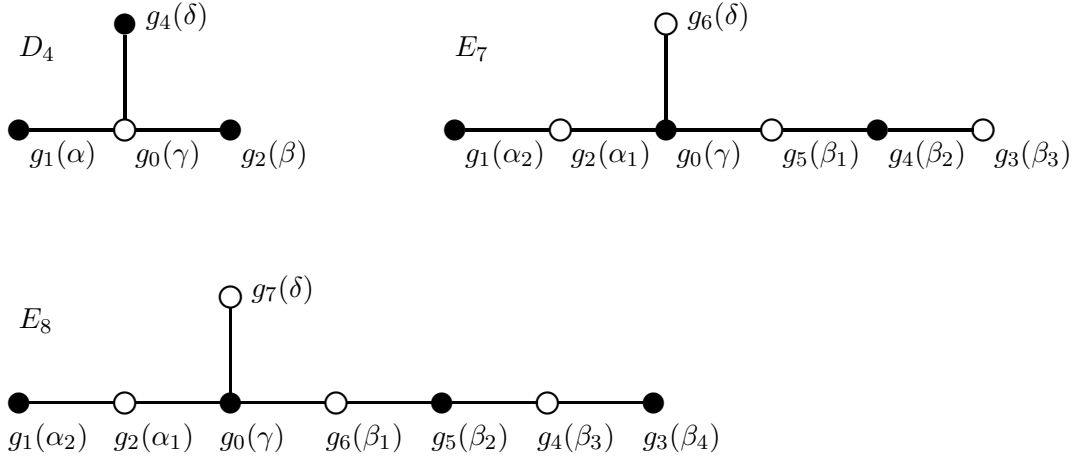
\begin{picture}(40,15)(-1,-1)
\thicklines
\put(-1,7){$D_4$}
 \drawline(9,1)(9,9)
 \multiputlist(0,0)(10,0){\circle*{2},\circle{2},\circle*{2}}
 \multiputlist(10,10)(0,10){\circle*{2}}
 \drawline(0,0)(8,0)
\drawline(10,0)(19,0) \put(0,-3){$g_1(\alpha)$}
\put(10,-3){$g_0(\gamma)$} \put(20,-3){$g_2(\beta)$}
\put(11,10){$g_4(\delta)$}
\end{picture}
 \vspace{5mm}
  \setlength{\unitlength}{4pt}
\begin{picture}(62,15)(-1,-1)
\thicklines
\put(-1,7){$E_7$}
 \drawline(19,0)(19,9)
\multiputlist(0,0)(10,0){\circle*{2},\circle{2},\circle*{2},\circle{2},\circle*{2},\circle{2}}
 \multiputlist(20,10)(0,10){\circle{2}}
 \drawline(0,0)(8,0)
\drawline(10,0)(19,0)
\drawline(20,0)(28,0)\drawline(30,0)(38,0)\drawline(40,0)(48,0)
\put(0,-3){$g_1(\alpha_2)$} \put(10,-3){$g_2(\alpha_1)$}
\put(20,-3){$g_0(\gamma)$} \put(30,-3){$g_5(\beta_1)$}
\put(40,-3){$g_4(\beta_2)$}
\put(50,-3){$g_3(\beta_3)$}\put(21,10){$g_6(\delta)$}
\end{picture}
\vspace{5mm}

 \setlength{\unitlength}{4pt}
\begin{picture}(62,15)(-1,-1)
\thicklines
\put(-1,7){$E_8$}
 \drawline(19,0)(19,9)
\multiputlist(0,0)(10,0){\circle*{2},\circle{2},\circle*{2},\circle{2},\circle*{2},\circle{2},\circle*{2}}
 \multiputlist(20,10)(0,10){\circle{2}}
 \drawline(0,0)(8,0)
\drawline(10,0)(19,0)
\drawline(20,0)(28,0)\drawline(30,0)(38,0)\drawline(40,0)(48,0)\drawline(50,0)(58,0)
\put(-2,-4){$g_1(\alpha_2)$} \put(8,-4){$g_2(\alpha_1)$}
\put(18,-4){$g_0(\gamma)$} \put(28,-4){$g_6(\beta_1)$}
\put(38,-4){$g_5(\beta_2)$}
\put(48,-4){$g_4(\beta_3)$}\put(58,-4){$g_3(\beta_4)$}\put(21,10){$g_7(\delta)$}
\end{picture}
\vspace{5mm}

This correspondence is important and henceforth we will designate the coordinates of the generalized dimension vectors by the vertices of the graph   $G$, $\Phi^{-1}(d)=(n_{g_0};n_{g_1},
n_{g_2},\ldots, n_{g_{k+l+m}})$.
 The following theorem holds. 
\begin{theorem}
The algebra  $\mathcal{P}_{\alpha,\beta,\delta,\gamma}$ associated with   $D_4$ has an irreducible non-degenerate *-representation only in the dimension    $(2;1,1,1)$. This representation is unique and exists  if and only if  the following conditions are  satisfied:
 \begin{gather*} \delta+\beta>\alpha, \alpha+\beta>\delta,\delta+\alpha>\beta,\\
2 \gamma=\alpha+\beta+\delta.
\end{gather*}
\end{theorem}

All irreducible non-degenerate locally-scalar representations of the graph  $E_7$ are obtained from the simplest representation   $\Pi_{g_0}$:  in 6 steps we obtain a non-degenerate  representation of dimension 
 $\overset{\circ}{c}_6(d_{g_0})=$ $(4; 1,2; 1,2,3; 2)$, and in 7 steps --- of dimension  $\overset{\circ}{c}_7(d_{g_0})=$ $(4; 1,3;
1,2,3; 2)$, in 8 steps  of dimension
$\overset{\circ}{c}_8(d_{g_0})=$ $(4; 2,3; 1,2,3; 2)$ (these are dimensions   $d$ of the  graph representations. Using them the  corresponding dimensions of the algebra representations can be computed as 
  $\Phi^{-1}(d)$). The following theorem holds. 

\begin{theorem}
The algebra  $\mathcal{P}_{\alpha,\beta,\delta,\gamma}$ associated with the  graph  $E_7$ has an irreducible non-degenerate  $*$-representation only in the following cases. 
\begin{enumerate}
\item   In dimension  $(4; 1,1; 1,1,1; 2)$, provided that
\begin{gather*}
\gamma >\delta +{{\beta }_3},\delta +{{\alpha }_1}+{{\beta }_1}> 2
\gamma ,\delta +{{\beta }_2}>\gamma , \\
 2 \gamma >\delta+{{\alpha
}_2}+{{\beta }_1}, \gamma
>{{\beta }_1}, 2 \gamma
>\delta +{{\alpha }_1}+{{\beta }_2},
\end{gather*}
 and $4\gamma=\alpha_1+\alpha_2+\beta_1+\beta_2+\beta_3+2 \delta$.
\item In dimension $(4; 1,2; 1,1,1; 2)$, provided that
\begin{gather*}
{{\alpha }_1}+ {{\beta }_2}>\gamma , 2  \gamma >{{\alpha
}_1}+{{\beta }_1}+\delta ,\gamma >
 {{\alpha }_1}+ {{\beta}_3},
 2 \gamma >\delta +{{\alpha }_1}+{{\beta }_2}+{{\beta }_3},\\
\delta +{{\alpha }_1}>\gamma ,
 \delta +{{\alpha}_1}+{{\beta }_1}+{{\beta }_3}>2 \gamma,
\end{gather*}
$4\gamma=2\alpha_1+\alpha_2+\beta_1+\beta_2+\beta_3+2 \delta$.
\item In dimension $(4; 2,1; 1,1,1; 2)$, provided that
\begin{gather*}
\gamma >{{\alpha }_2}+{{\beta }_2},\delta +{{\alpha }_2}+{{\beta
}_3}>\gamma ,{{\alpha }_2}+{{\beta }_1}>\gamma ,
 \delta +{{\alpha }_2}+{{\beta }_1}+{{\beta }_2}>2 \gamma,\\
\gamma >\delta +{{\alpha }_2},
2 \gamma >\delta +{{\alpha }_2}+{{\beta }_1}+{{\beta }_3},
\end{gather*}
$4\gamma=\alpha_1+2 \alpha_2+\beta_1+\beta_2+\beta_3+2 \delta$.
\end{enumerate}
In each of the above cases, the representation is unique up to unitary equivalence.  
\end{theorem}

Non-degenerate irreducible representations of the graph 
 $E_8$ are obtained from the simplest ones   $\Pi_{g_0}$ and   $\Pi_{g_6}$.  In  $k$ steps, where  $k=9,10, \ldots, 14$, starting from  $\Pi_{g_0}$ we obtain a non-degenerate representation of dimension    $\overset{\circ}{c}_k(d_{g_0})$.
 Analogously, in $k$ steps, where $k=11, 12, 13, 14$, using   $\Pi_{g_6}$
we obtain a non-degenerate representation of dimension
$\overset{\circ}{c}_k(d_{g_6})$, where  $d_{g_6}(g_6)=1$ and 
$d_{g_6}(g_i)=0$, if  $g_i\not= g_6$.

\begin{theorem}\label{ineqe8}
The algebra $\mathcal{P}_{\alpha,\beta,\delta,\gamma}$ associated with graph   $E_8$ has an irreducible non-degenerate *-representation only in the following cases. 
\begin{enumerate}
\item   In dimension $(5;2,1;1,1,1,1;2)$, provided that
\begin{gather*}
\gamma >{{\alpha }_2}+{{\beta }_3}, 
 \delta +{{\alpha }_2}+{{\beta }_1}+{{\beta }_4}>2\gamma,  
{{\alpha }_2}+{{\beta }_2}>\gamma,\\ 
 \delta +{{\alpha }_2}+{{\beta }_2}+{{\beta }_3}>2\gamma,   
2\gamma >\delta +{{\alpha }_2}+{{\beta }_1},\gamma >{{\beta }_1}, \\
2\gamma >\delta   
+{{\alpha }_2}+{{\beta }_2}+{{\beta }_4},\\
\gamma =
 \frac{1}{5} (2 \delta +{{\alpha }_1}+2 {{\alpha }_2}+{{\beta }_1}+{{\beta }_2}+{{\beta }_3}+{{\beta }_4}).
\end{gather*}
\item In dimension $(5;2,2;1,1,1,1;3)$, provided that
\begin{gather*}
\delta +{{\alpha }_2}+{{\beta }_3}>\gamma,   
 5 \gamma +{{\beta }_1}+{{\beta }_3}+{{\beta }_4}>  
 2 \delta +3 {{\alpha }_1}+3 {{\alpha }_2}+4 {{\beta }_2},\\  
 \gamma >\delta +{{\alpha }_2}+{{\beta }_4},   
 4 \delta +{{\alpha }_1}+{{\alpha }_2}+3 {{\beta }_1}+3 {{\beta }_2}>  
5 \gamma +2 {{\beta }_3}+2 {{\beta }_4},\\  
 2 \delta +3 {{\alpha }_1}+3 {{\alpha }_2}+4 {{\beta }_1}>  
 5 \gamma +{{\beta }_2}+{{\beta }_3}+{{\beta }_4},\delta +{{\alpha }_1}>\gamma,\\
 2 \delta +3 {{\alpha }_1}+3 {{\alpha }_2}+4 {{\beta }_2}+4 {{\beta }_4}>  
 5 \gamma +{{\beta }_1}+{{\beta }_3}, \\
\gamma = \frac{1}{5} (3 \delta +2 {{\alpha }_1}+2 {{\alpha
}_2}+{{\beta }_1}+{{\beta }_2}+
  {{\beta }_3}+{{\beta }_4}).
\end{gather*}
\item In dimension $(6;2,2;1,1,1,1;3)$, provided that
\begin{gather*}
2 \gamma >\delta +{{\alpha }_1}+{{\beta }_2},   
 \delta +{{\alpha }_1}+{{\alpha }_2}+{{\beta }_3}>2 \gamma,   
 \delta +{{\alpha }_1}+{{\beta }_1}>2 \gamma,\\   
 3 \gamma >2 \delta +{{\alpha }_1}+{{\alpha }_2}+{{\beta }_3}+{{\beta }_4},   
 2 \gamma >\delta +{{\alpha }_1}+{{\alpha }_2}+{{\beta }_4},   
 \gamma >\delta +{{\alpha }_2}, \\  
 2 \delta +{{\alpha }_1}+{{\alpha }_2}+{{\beta }_2}+{{\beta }_4}>3 \gamma, \\
\gamma =   \frac{1}{6} (3 \delta +2 {{\alpha }_1}+2 {{\alpha
}_2}+{{\beta }_1}+{{\beta }_2}+
 {{\beta }_3}+{{\beta }_4}).
\end{gather*}
\item In dimension $(6;2,2;1,1,1,2;3)$, provided that
\begin{gather*}
\delta +{{\alpha }_1}+{{\beta }_1}+{{\beta }_4}>2 \gamma ,   
 2 \gamma >\delta +{{\alpha }_2}+{{\beta }_1}+{{\beta }_3},   
 2 \gamma >\delta +{{\alpha }_1}+{{\beta }_1}, \\  
 3 \gamma >2 \delta +{{\alpha }_1}+{{\alpha }_2}+{{\beta }_1}+{{\beta }_4},  
 \delta +{{\alpha }_2}+{{\beta }_1}+{{\beta }_2}>2 \gamma,   
 {{\alpha }_2}+{{\beta }_1}>\gamma,\\   
 2 \delta +{{\alpha }_1}+{{\alpha }_2}+{{\beta }_1}+{{\beta }_3}>3 \gamma, \\
\gamma =  \frac{1}{6} (3 \delta +2 {{\alpha }_1}+2 {{\alpha }_2}+2
{{\beta }_1}+{{\beta }_2}+
 {{\beta }_3}+{{\beta }_4}).
\end{gather*}
\item In dimension $(6;2,2;2,1,1,1;3)$, provided that
\begin{gather*}
2 \gamma >\delta +{{\alpha }_2}+{{\beta }_1}+{{\beta }_4},   
 \delta +{{\alpha }_1}+{{\beta }_2}+{{\beta }_4}>2 \gamma,   
 \delta +{{\alpha }_2}+{{\beta }_4}>\gamma, \\    
 2 \delta +{{\alpha }_1}+{{\alpha }_2}+{{\beta }_1}+{{\beta }_4}>3 \gamma ,   
 2 \gamma >\delta +{{\alpha }_1}+{{\beta }_3}+{{\beta }_4},   
 \gamma >{{\alpha }_1}+{{\beta }_4}, \\  
 3 \gamma >2 \delta +{{\alpha }_1}+{{\alpha }_2}+{{\beta }_2}+{{\beta }_4}, \\
\gamma =    \frac{1}{6} (3 \delta +2 {{\alpha }_1}+2 {{\alpha
}_2}+{{\beta }_1}+{{\beta }_2}+
 {{\beta }_3}+2 {{\beta }_4}).
\end{gather*}
\item In dimension $(6;2,2;1,1,2,1;3)$, provided that
\begin{gather*}
\delta +{{\alpha }_1}+{{\alpha }_2}+{{\beta }_2}>2 \gamma ,   
 2 \gamma >\delta +{{\alpha }_1}+{{\beta }_2}+{{\beta }_4},   
 \gamma >{{\alpha }_2}+{{\beta }_2}, \\ 
 3 \gamma >\delta +{{\alpha }_1}+{{\alpha }_2}+{{\beta }_1}+{{\beta }_2},   
 \delta +{{\alpha }_1}+{{\beta }_2}+{{\beta }_3}>2 \gamma ,   
 \delta +{{\beta }_2}>\gamma ,\\   
 3 \gamma >2 \delta +{{\alpha }_1}+{{\alpha }_2}+{{\beta }_2}+{{\beta }_3}, \\
 \gamma =  \frac{1}{6} (3 \delta +2 {{\alpha }_1}+2 {{\alpha }_2}+{{\beta }_1}+2 {{\beta }_2}+
  {{\beta }_3}+{{\beta }_4}).
\end{gather*}
\item In dimension $(6;2,2;1,2,1,1;3)$, provided that
\begin{gather*}
2 \gamma >\delta +{{\alpha }_1}+{{\alpha }_2}+{{\beta }_3},   
 \delta +{{\alpha }_2}+{{\beta }_1}+{{\beta }_3}>2 \gamma,   
 {{\alpha }_1}+{{\beta }_3}>\gamma, \\    
 \delta +{{\alpha }_1}+{{\alpha }_2}+{{\beta }_3}+{{\beta }_4}>2 \gamma ,   
 2 \gamma >\delta +{{\alpha }_2}+{{\beta }_2}+{{\beta }_3},   
 \gamma >\delta +{{\beta }_3}, \\ 
 2 \delta +{{\alpha }_1}+{{\alpha }_2}+{{\beta }_2}+{{\beta }_3}>3 \gamma, \\
\gamma =  \frac{1}{6} (3 \delta +2 {{\alpha }_1}+2 {{\alpha
}_2}+{{\beta }_1}+{{\beta }_2}+ 2 {{\beta }_3}+{{\beta }_4}).
\end{gather*}
\item In dimension $(5;1,2;1,1,1,1;2)$, provided that
\begin{gather*}
5 \gamma +{{\alpha }_2}+{{\beta }_1}+{{\beta }_3}+{{\beta }_4}>  
 3 \delta +3 {{\alpha }_1}+4 {{\beta }_2},{{\alpha }_1}+{{\beta }_3}>\gamma ,   \gamma >{{\beta }_1},\\ 
5 \gamma +2 {{\beta }_1}>  
 \delta +{{\alpha }_1}+3 {{\alpha }_2}+3 {{\beta }_2}+3 {{\beta }_3}+3 {{\beta }_4}, \\ 
 5 \gamma +{{\alpha }_2}+{{\beta }_1}+{{\beta }_2}>  
 3 \delta +3 {{\alpha }_1}+4 {{\beta }_3}+4 {{\beta }_4},  
 \gamma >{{\alpha }_1}+{{\beta }_4}, \\ 
 5 \gamma +2 {{\beta }_2}+2 {{\beta }_4}>  
 \delta +{{\alpha }_1}+3 {{\alpha }_2}+3 {{\beta }_1}+3 {{\beta }_3}, \\
\gamma =   \frac{1}{5} (2 \delta +2 {{\alpha }_1}+{{\alpha
}_2}+{{\beta }_1}+{{\beta }_2}+{{\beta }_3}+{{\beta }_4}).
\end{gather*}
\item In dimension $(5;1,2;1,1,1,1;3)$, provided that
\begin{gather*}
5 \gamma +2 {{\beta }_1}+2 {{\beta }_4}>  
 4 \delta +{{\alpha }_1}+3 {{\alpha }_2}+3 {{\beta }_2}+3 {{\beta }_3},  
\gamma >\delta +{{\beta }_3},\\ \delta +{{\alpha }_1}>\gamma,  
 5 \gamma +{{\alpha }_2}+{{\beta }_2}+{{\beta }_3}+{{\beta }_4}>  
 2 \delta +3 {{\alpha }_1}+4 {{\beta }_1},\\  
 5 \gamma +2 {{\beta }_2}+2 {{\beta }_3}>  
 4 \delta +{{\alpha }_1}+3 {{\alpha }_2}+3 {{\beta }_1}+3 {{\beta }_4},  
 \delta +{{\beta }_2}>\gamma, \\  
 5 \gamma +{{\alpha }_2}+{{\beta }_1}+{{\beta }_3}>  
 2 \delta +3 {{\alpha }_1}+4 {{\beta }_2}+4 {{\beta }_4}, \\
\gamma =   \frac{1}{5} (3 \delta +2 {{\alpha }_1}+{{\alpha
}_2}+{{\beta }_1}+{{\beta }_2}+{{\beta }_3}+{{\beta }_4}).
\end{gather*}
\item In dimension $(5;2,1;1,1,1,1;3)$, provided that
\begin{gather*}
5 \gamma +{{\alpha }_1}+{{\beta }_2}+{{\beta }_3}>  
2 \delta +3 {{\alpha }_2}+4 {{\beta }_1}+4 {{\beta }_4}, \delta +{{\beta }_2}>\gamma, \\ 
 \gamma >\delta +{{\alpha }_2},\delta +{{\alpha }_2}+{{\beta }_4}>\gamma,    
 5 \gamma +{{\alpha }_1}+{{\beta }_1}+{{\beta }_4}>  
 2 \delta +3 {{\alpha }_2}+4 {{\beta }_2}+4 {{\beta }_3},\\  \gamma >\delta +{{\beta }_3},   
 5 \gamma +2 {{\beta }_1}+2 {{\beta }_3}>  
 4 \delta +3 {{\alpha }_1}+{{\alpha }_2}+3 {{\beta }_2}+3 {{\beta }_4}, \\
\gamma =  \frac{1}{5} (3 \delta +{{\alpha }_1}+2 {{\alpha
}_2}+{{\beta }_1}+{{\beta }_2}+{{\beta }_3}+{{\beta }_4}).
\end{gather*}
\item In dimension $(5;2,2;1,1,1,1;2)$, provided that
\begin{gather*}
5 \gamma +2 {{\beta }_2}>\delta +{{\alpha }_1}+{{\alpha }_2}+3 {{\beta }_1}+  
 3 {{\beta }_3}+3 {{\beta }_4},\gamma >{{\alpha }_1}+{{\beta }_4}, \\ 
 {{\alpha }_2}+{{\beta }_1}>\gamma ,\gamma >{{\alpha }_2}+{{\beta }_2},  
 5 \gamma +2 {{\beta }_3}+2 {{\beta }_4}>  
 \delta +{{\alpha }_1}+{{\alpha }_2}+3 {{\beta }_1}+3 {{\beta }_2}, \\ 
 {{\alpha }_1}+{{\beta }_3}>\gamma,  
 5 \gamma +{{\beta }_1}+{{\beta }_2}+{{\beta }_4}>  
 3 \delta +3 {{\alpha }_1}+3 {{\alpha }_2}+4 {{\beta }_3}, \\
\gamma =   \frac{1}{5} (2 \delta +2 {{\alpha }_1}+2
{{\alpha}_2}+{{\beta }_1}+{{\beta }_2}+
   {{\beta }_3}+{{\beta }_4}).
\end{gather*}
\end{enumerate}
In each of the above cases the representation is unique up to unitary equivalence.
\end{theorem}

\noindent {\it Remark.}   Inequalities in Theorem~\ref{ineqe6}-- \ref{ineqe8} can be obtained alternatively using Horn's inequalities (see, for example~\cite{fulton}) and W. Crawley-Boevey solution of the  Deligne-Simpson problem (see~\cite{cwduke}).  More precisely, representations of   $\mathcal{P}_{\alpha,\beta,\delta,\gamma}$ (not necessarily $*$-preserving) in a fixed generalized dimension $d=(n;\ldots )$ are in one-to-one correspondence with  $n\times n$ matrices $A_1$, $A_2$, $A_3$  picked from fixed diagonal  conjugacy classes $C_1$, $C_2$, $C_3$. These classes are determined by $(\alpha, \beta, \gamma, \delta)$ and $d$. 
 Thus an application  of \cite[Theorem~1]{cwduke} to diagonal conjugacy classes gives a description of the  parameters $(\alpha, \beta, \gamma, \delta)$  for which  there is an irreducible  (not necessarily $*$-preserving) representation of algebra  $\mathcal{P}_{\alpha,\beta,\delta,\gamma}$ in terms of the root system of a certain Kac-Moody Lie algebra. In particular, it follows that there are only finitely many such representations in the case of a Dynkin graph and generalized dimensions $d$ are some simple roots (up to a certain linear transformation (\cite[Theorem~1]{cwduke}). Moreover, for a fixed generalized dimension $d$ such a representation $\pi$ is unique.  If we require that the parameters of the algebra satisfy Horn's inequalities (we must translate the parameters $\alpha, \beta, \gamma, \delta$ and the generalized dimension $d$ into eigenvalues and their multiplicities first) then $\pi$ must be unitarizable and thus $(\alpha, \beta, \gamma, \delta)\in W_{irr}$.  Thus combining Horn's inequalities and  Crawley-Boevey conditions of irreducibility for $d$  we obtain inequalities in Theorem~\ref{ineqe6}- \ref{ineqe8}. However, to derive our results we do not use neither Horn's inequalities nor Crawley-Boevey conditions. Moreover, we give an algorithm, the Coxeter functors for $*$-representations, to construct $*$-representations of  $\mathcal{P}_{\alpha,\beta,\delta,\gamma}$.

}

\end{document}